\documentclass{article}
\usepackage{hyperref}
\hypersetup{
    colorlinks=true,
    linkcolor=red,
    urlcolor=blue,
    citecolor=green}
\usepackage{arxiv}
\usepackage[english]{babel}
\usepackage[utf8]{inputenc} % allow utf-8 input
\usepackage[T1]{fontenc}    % use 8-bit T1 fonts
\usepackage{hyperref}       % hyperlinks
\usepackage{url}            % simple URL typesetting
\usepackage{booktabs}       % professional-quality tables
\usepackage{amsfonts}       % blackboard math symbols
\usepackage{nicefrac}       % compact symbols for 1/2, etc.
\usepackage{microtype}      % microtypography
\usepackage{lipsum}
\usepackage{amsmath}
\usepackage{amssymb}
\usepackage{pdfpages}
\usepackage[square,numbers]{natbib}
\usepackage[nottoc]{tocbibind}
\usepackage[T1]{fontenc}
\usepackage{natbib}
\title{Primal Conjecture in Matrix $E_\phi^9$.}

\author{
Paul Marrero \\
  \texttt{Paulqed0@protonmail.com} \\
\And
 Eduardo J. Acuña T.\\
  \texttt{eacuna1@uc.edu.ve} \\}

\date{April 2021}

\begin{document}

\maketitle
\begin{abstract}
\begin{center}
    \emph{"There is a certain enigma in everything that seems obvious."}
\end{center}

   In this paper we propose a conjecture about integer solutions to any equations, based on Primal algebra specifically this conjecture is a corollary of the Acuña´s Theorem in that article.\\
   Also some problems are proposed which, if the conjecture is correct, could be solved.
   
\end{abstract}

\section{Introduction}
The approach of this article is from the mathematical philosophical position of constructivism, that is to say that we create an arithmetic system with our rules and we observe how the elements behave and what relation they have with the real mathematics.\\
With this in mind we introduce Primordial Algebra as a system of analysis aimed at integers, especially operations between integers $\pm \mathbb{Z}$.\\
This algebra analyzes the numbers under a unitary module which is the 9 in base 10, in this way the remainders of each positive or negative integer generate groups congruent to each other, although these properties have been studied in the past it had not been possible to operate between positive and negative, the primordial algebra allows to operate between both types without falling into incongruities.
\begin{align*}
    \frac{-Z}{-9Z}~ \bigcup~ \frac{Z}{9Z}~ /~ Z \in~ \mathbb{Z}^+,~-Z ~\in~ \mathbb{Z}^-
\end{align*}
\subsection{Notation}
Before beginning to understand how primordial algebra works, it is necessary to understand the notation created for it.\\

\begin{itemize}
    \item Let $v$ be any integer number
    \item Let $\Phi$ be a function of the digital root of any $\mathbb{Z}$
    \item Let $\Phi^Z$ be any of the class of any $\mathbb{Z}$, we replace $"Z"$ to denote the value of the class.
    \item Let $Z_\phi$ be any Primal Number, which is the representative value of the class and can be operated with.
\end{itemize}
\subsection{Reduction of integers}
The First idea that we use for this matrix was the property of reduction of integer numbers by the Digital Sum, all this is mention in the article (2021)  \cite{tarazona2021matrize_phi},
This is congruent whit the Method of Digital Sum  Propose by Jhon Eider (2015) \cite{lopez2015suma}, and we propose the the next Equation to found any Positive Digital Root based on Knuth (1972) \cite{knuth2005art}.\\
Also this paper is gonna drop a resume of primal algebra, with the propose of make more easy to understand the conjecture.

The Novem Digital Root:  

\begin{equation}
\label{equ1}
    \Phi_9=v-9\left(\left\lfloor\frac{v-1}{9}\right\rfloor\right)
\end{equation}

in the Matrix $E_\phi^9$ the equation for the Digital Sum is:

\begin{equation}
\label{equ2}
    \Phi_j=v-9(i-1)
\end{equation}

Both of this equations are equivalent, so we conclude when the digital root is in module 9, is equal to the Digital Sum.

\section{Theorem Marrero Criterion}
The Marrero's criterion, (2021) \cite{tarazona2021matrize_phi}, also mentioned in the pre-print \cite{https://doi.org/10.48550/arxiv.2103.17037} as the $K$-interactive of the digital Sum.
This is a criterion who tell when to stop to the Digital Sum $\Phi$, we denoted the Criterion like "$M$" this is the mark to any function of Digital Sum that meet the condition of the criterion $\phi^M$.
\subsection{Theorem M}
Let $v$ are integer, then there is a unique integer $d \in \hspace{0.1cm} D$ such that:
\begin{equation}
\label{equ3}
    \phi^k(v_k)- d = 0,
\end{equation}
where $D$ of the set of the points fixed of $v$.
\subsection{Proof:}

\begin{enumerate}
    \item (Existence). Supposed $b \in \mathbb Z$ and $k$ be a integer fixed. Consider the set $ S \left\{ b = \phi^k(v_k)-b=0 \right\}$ where $v_k \in \mathbb Z $, $S$ is not empty, because $1 \in S $ with $k=1$.
Then $S \subseteq \mathbb Z$ and let $d$ being the smallest element of $S, d \in S$, it must have the form $\phi^k(v_k)-d=0$, therefore $\phi^k(v_k)=d.$

    \vspace{0.4cm}
    
    \item (Uniqueness). Let $d_1$ and $d_2 \in D$ such that $d_1 \neq d_2$. supposed that $\phi^{k_1}(v_k)- d_1 = 0$ and $\phi^{k_2}(v_k)- d_2 =0$ (using Definition 2.1), then
$\phi^{k_1}(v_k)- d_1 = \phi^{k_2}(v_k)- d_2,$ but $\phi^{k_1}(v_k)-\phi^{k_2}(v_k)= d_1-d_2$ and $\phi^{k_1}(v_k)-\phi^{k_2}(v_k)=0$. therefore $d_1=d_2$, this constitutes a contradiction.  $\square$
\end{enumerate}

\vspace{1cm}
\section{Elements of the Matrix}
In this section we give a brief introduction to the matrix $E_\phi^9$, it is formed by the union of two positive and negative sub-matrices together with its own classification system, we adopt the following convention: Primal number $Z_\phi$ ,Primal class $\Phi^Z$, integer numbers $v$, rows $(i)$, columns $(j)$. See the following figures for more details.

\begin{center}

    $\begin{bmatrix}
1_\phi & 2_\phi  & 3_\phi & 4_\phi & 5_\phi& 6_\phi & 7_\phi & 8_\phi & 9_\phi\\
1 & 2 &  3 & 4 & 5 & 6 & 7 & 8 & 9 \\ 
10 & 11 & 12 & 13 & 14 & 15 & 16 & 17 & 18 \\
19 & 20 & 21 & 22 & 23 & 24 & 25 & 26 & 27\\
28 & 29 & 30 & 31 & 32 & 33 & 34 & 35 & 36 \\
\vdots & \vdots & \vdots & \vdots & \vdots & \vdots & \vdots & \vdots & \vdots \\
\end{bmatrix}$

Table 1: \emph{Positive sub-matrix} $[V^+]$ 
\end{center}

\begin{center}
    $\begin{bmatrix}
-1_\phi & -2_\phi  & -3_\phi & -4_\phi & -5_\phi& -6_\phi & -7_\phi & -8_\phi & -9_\phi\\
-1 & -2 &  -3 & -4 & -5 & -6 & -7 & -8 & -9 \\ 
-10 & -11 & -12 & -13 & -14 & -15 & -16 & -17 & -18 \\
-19 & -20 & -21 & -22 & -23 & -24 & -25 & -26 & -27\\
-28 & -29 & -30 & -31 & -32 & -33 & -34 & -35 & -36 \\
\vdots & \vdots & \vdots & \vdots & \vdots & \vdots & \vdots & \vdots & \vdots \\
\end{bmatrix}$

Table 2: \emph{Negative sub-matrix} $[V^-]$
\end{center}

\vspace{1cm}
\subsection{Primal Numbers}

The primal number is the result of applying the digital sum function $\Phi$ and we will denote it by $Z_\phi$. We denote the set of positive and negative primal numbers by 
 
 \begin{center}
     $\begin{matrix}
     D=\left\{{{d|-9_\phi \leq d \leq 9_\phi  |}}\right\}
 \end{matrix}$
 \end{center}
    
Note that each value of $D$ represents a relationship between and its classes, that is, for every integer $v$ that belongs to column $j$, there exists a primal number $Z_\phi$ such that $Z_\phi$ is congruent with $j$ module 9.

\subsubsection{Example}
\begin{center}
    \begin{enumerate}
        \item \hspace{5cm} $9_\phi = 72$
        \item \hspace{5cm} $8_\phi = 179$
        \item \hspace{5cm} $4_\phi = 9211$ \\
\hspace{5.6cm} \vdots
        \item \hspace{4.8cm} $-5_\phi = -14$
        \item \hspace{4.8cm} $-9_\phi = -927$
    \end{enumerate}
\end{center}

 \subsection{The Null Primal class}
  The Null class $\emptyset$, we will denote it by 0, is the one that is constituted by a column of zeros, and establishes a numerical border between the positive and negative classes of the matrix $E_\phi^9$.
  \vspace{1cm}
  \begin{center}
       $\begin{bmatrix}
-9_\phi & -8_\phi  & \ldots & -1_\phi & \emptyset & 1_\phi & 2_\phi & \ldots & 9_\phi\\
-9 & -8 & \ldots & -1 & 0 &  1 & 2 & \ldots & 9 \\ 
-18 & -17 & \ldots & -10 & 0 & 10 & 11 & \ldots & 18 \\
-27 & -26 & \ldots & -19 & 0 & 19 & 20 & \ldots & 27\\
-36 & -35 & \ldots & -28 & 0 & 28 & 29 & \ldots & 36 \\
-45 & -44 & \ldots & -37 & 0 & 37 & 38 & \ldots & 45 \\
-54 & -53 & \ldots & -46 & 0 & 46 & 47 & \ldots & 54 \\
-63 & -62 & \ldots & -55 & 0 & 55 & 56 & \ldots & 63 \\
-72 & -71 & \ldots & -64 & 0 & 64 & 65 & \ldots & 72 \\
-81 & -80 & \ldots & -73 & 0 & 73 & 72 & \ldots & 81 \\
-90 & -89 & \ldots & -82 & 0 & 82 & 83 & \ldots & 90 \\
-99 & -98 & \ldots & -91 & 0 & 91 & 92 & \ldots & 99 \\
\vdots & \vdots & \vdots & \vdots & \vdots & \vdots & \vdots & \vdots & \vdots \\
\end{bmatrix}$  

Table 3: Matrix $E_\phi^9$
  \end{center}

\subsection{Primal Class}
    The primal class that we denote by $\Phi^Z$ represents and are sets of congruence classes modulo 9, that is:

    \begin{center}
        $\Phi^Z= \left\{v \in \mathbb Z \mid Z_\phi \equiv v \hspace{0,1cm}(mod~9)  \right\} $
    \end{center}

\subsubsection{Example:} 

 \begin{enumerate}
     \item \begin{center}
         $\Phi^1 \hspace{0.2cm}= \left\{1,10,19,28,...,91...\right\}$
     \end{center}
     \item \begin{center}
         $\Phi^2 \hspace{0.2cm}= \left\{2,11,20,29,...,92...\right\}$
     \end{center}
     \item \begin{center}
            $\Phi^9 \hspace{0.2cm}= \left\{9,18,27,36,...,99...\right\}$
     \end{center}
     \item \begin{center}
         $\Phi^{-1} = \left\{-1,-10,...,-199...\right\}$
     \end{center}
     \item \begin{center}
     $\Phi^{-9} = \left\{-9,-18,...,-999...\right\}$
     \end{center}
 \end{enumerate}

 The Primal class also are class of congruence with respect to the module 9, Vinogradov (1977) \cite{vinogradov1977fundamentos}, in class $9^+_-$ the residue is 0, but we take the module as the main value of the class, either positive or negative.
\vspace{1.5cm}
 
 \section{Method of calculation of the positive columns}
 
 Recall that the columns of the matrix $E_\phi^9$, we denote it by $j$, it is worth mentioning that said Matrix $E_\phi^9$,consists of 18 columns, nine are for the positive part $[V^+]$ and nine for the negative part $[V^-]$ along with its respective primal number that we write at the beginning of each column in the matrix $E_\phi^9$.

There is many methods to calculate the columns by definition we chose the Digital Sum $\Phi$, The novem digital$\Phi_9$\ref{equ1} or The matrix form $\Phi_j$\ref{equ2}.
All mentioned before are valid to the positive Part $[V^+]$.

\section{ Method of calculation of the Negative Columns}
As the same way we define the positive columns, the negative part of the matrix contains all the negative integers $\mathbb Z^-$, this calculate is possible thanks to the negative novem digital root.
\subsection{Definition:}
Let $v$ be a negative integer, we will write it as $-v$. We define the \textbf{negative digital root novem}, which we denote by $\mho$ and is as follows:

\begin{equation}
\label{equ4}
    \mho_{-9}(v) = ((v+1)~mod\hspace{0.1cm} (-9)) -1
\end{equation}

\subsection{Definition 2:}
We define the \textbf{negative novem floor function},\cite{paul_francisco_marrero_romero_2021_5598339}  for a negative integer, we will denote it by  and it is as follows:

\begin{equation}
\label{equ5}
   \hat{\mho}_{-9}(v)= v-(-9) {\left\lfloor{\frac{v+1}{-9}}\right\rfloor}
\end{equation}

For the Negative Digital Sum in the matrix $\Phi_{-j}$ we propose the next equation:

\begin{equation}
\label{equ6}
    \Phi_{-j}=v+9(i-1)
\end{equation}
 Where equation \ref{equ5} and \ref{equ6} are congruent with each other.
 \vspace{2cm}
 
 \section{Method of Calculation of the rows of the Matrix.}
 Let $i$ be the set of rows of the Matrix $E_\phi^9$, the rows of the Matrix are always positive for any number inside the matrix, no matters if is negative, but the method to calculate their $i$ is conditional to their sing of the number.
 
 To calculate the rows of the Matrix $E_\phi^9$ we proceed to consider two cases:
 
 \begin{itemize}
     \item case 1: Row a positive integer.
     
      \vspace{0.3cm}
     \begin{enumerate}
         \item[]
         \begin{equation}
         \label{equ7}
      i_{v}=\left\lbrace\begin{array}{c}\left\lfloor\frac{v}{9}+1\right\rfloor \quad~For~a~v\neq9c \\ 
      \vspace{0.2cm} \\ \frac{v}{9} \hspace{1.2cm}~For~a~v =9c \end{array}\right.
  \end{equation}
     \end{enumerate}
       \item case 2: Row a negative integer.
       
       \vspace{0.3cm}
       \begin{enumerate}
           \item[]
           \begin{equation}
           \label{equ8}
      i_{-v}=\left\lbrace\begin{array}{c}\left\lfloor-\frac{v}{9}+1\right\rfloor \quad~For~a~v\neq9c \\
      \vspace{0.2cm} \\ -\frac{v}{9} \hspace{1.2cm}~For~a~v =9c \end{array}\right.
  \end{equation}
       \end{enumerate}
 \end{itemize}
  \vspace{0.2cm}
  
 Where $9c$ are all the multiple of 9.
 
 Other ways to calculate the $i$ is whit this two equations:
 \begin{equation}
 \label{equ9}
     i_v=\frac{v-\phi_v}{9}+1
 \end{equation}
    
\begin{equation}
\label{eque10}
      i_{-v}=\frac{v-\phi_{-v}}{9}-1
\end{equation}
  
  \section{Decompress Primal numbers}
  One of the advantage of putting the integers numbers in the Matrix $E_\phi^9$ and turning into Primal numbers is that every number no matter the size can be compress into one number, but whit the matrix equations we can decompress those Primal numbers and get back the integer number.
  This is the fundamental base for the creation of a algebra of from the matrix, because the equations to decompress are getting isolating the formulas \ref{equ2} and \ref{equ6}:
  \vspace{0.2cm}

  Isolating equation \ref{equ2} to get the positive number we get:
  \vspace{0.2cm}
     \begin{equation}
     \label{equ11}
         v=\Phi_j+9(i-1)
     \end{equation}
\vspace{0.2cm}
     Isolating equation \ref{equ6} to get the negative number we get:
\vspace{0.2cm}
     \begin{equation}
     \label{equ12}
         -v=\Phi_{-j}-9(i-1)
     \end{equation}

\subsection{Notation for reduced numbers in the matrix.}
As we know a integer number in the Matrix is a set of three elements, his Primal Number $Z_\phi$, his position in $j$ and $i$, whit this we can represent a integer inside the matrix like this:
\begin{align*}
    v=Z_i
\end{align*}
Where $i$ is the row and $Z$ is the representative number, we don't use $\phi$ because would be redundant such that $\phi \equiv Z \equiv j$ are numerically the same, so we prefer pick only one to represent all.
\subsubsection{Example:} Now some examples of how represent numbers in the matrix.

\begin{enumerate}

    \item \begin{center}
    $v=178,\rightarrow Z_i = 7_{20}$
    \end{center}
    \item \begin{center}
   $v=157,\rightarrow Z_i = 4_{18}$
    \end{center}
    \item \begin{center}
        $v=293,\rightarrow Z_i = 5_{33}$
    \end{center}
    \item \begin{center}
        $v=359,\rightarrow Z_i = 8_{40}$
    \end{center}
    \item \begin{center}
        $v=999,\rightarrow Z_i = 9_{111}$
    \end{center}
\end{enumerate}
If we go to the Matrix and locate the $Z$ in the columns and the $i$ in the rows the interception would be the correspondent number.
Now whit the decompress equation we can go from $Z_i$ to $v$ using the example 4. shown previously.
\vspace{0.2cm}
\begin{align*}
    \begin{matrix}
        8_{40}=8+9(40-1) \rightarrow 8_{40}=8+9(39) \\
        8_{40}=8+351 \rightarrow 8_{40}=359\\
    \end{matrix}
\end{align*}
\vspace{0.2cm}
 \section{The Acuña's Theorem}
 \label{acuña}
 The Acuña's Theorem \cite{tarazona2021matrize_phi}, is the core theorem that expose the principles of the algebra primordial, using the properties of the matrix to arithmetic between classes.
 \subsection{Definition:}
 
  For all operation that we denote as $(\odot=+,-,/,*)$ between two or more Primal Classes, the result will be another Primal class, and this is congruent with the numbers inside the class.
 
 \begin{align*}
     \forall \hspace{0.2cm} Z_{\phi a} \odot Z_{\phi b} \equiv  v_a \odot \hspace{0.1cm} v_b~(mod~9)
 \end{align*}
 
 This is a extension of of class of residues by Vinogradov (1977) \cite{vinogradov1977fundamentos}, but this work in both positive and negative integers.
 
 \subsection{Proof:}
 If:
 \begin{align*}
   \begin{matrix}
    \vspace{0.2cm}
       v_a, Z_{\phi a} \in \phi^a \rightarrow v_a \equiv  Z_{\phi a}~(mod\hspace{0.1cm} 9)\\
       \vspace{0.2cm}
       v_b, Z_{\phi b} \in \phi^b \rightarrow v_b \equiv  Z_{\phi b}~(mod\hspace{0.1cm} \hspace{0.1cm} 9)\\
       v_c, Z_{\phi c} \in \phi^c \rightarrow v_c \equiv  Z_{\phi c}~(mod\hspace{0.1cm} 9)
   \end{matrix}
 \end{align*}
 Then by definition:
 
 \begin{align*}
     \begin{matrix}
      (\phi \odot \phi)^{a+b=}=\left\{ x \in \mathbb Z \mid x \equiv Z_{\phi a} + Z_{\phi b} = Z_{\phi c}~ (mod\hspace{0.1cm} 9) \right\}
     \end{matrix}
 \end{align*}
 From where:
 
 \begin{align*}
     \begin{matrix}
     \vspace{0.2cm}
      x-(Z_{\phi a}+Z_{\phi b}=Z_{\phi c}) = 9t, \hspace{0.2cm} for\hspace{0.1cm} t \in \mathbb Z\\ 
      \vspace{0.2cm}
      x=9t+(Z_{\phi a}+Z_{\phi b}=Z_{\phi c})\\
      \vspace{0.2cm}
      x=9t+[(v_a - 9q)+ (v_b - 9r) = (v_c - 9s)] \hspace{0.2cm}
      for\hspace{0.2cm} q, r, s \in \mathbb Z \\
      \vspace{0.2cm}
      x=9t+[(v_a + v_b)-9(q+r)= v_c - 9s]\\
      \vspace{0.2cm}
      x=9t+[(v_a + v_b) - 9(q+r) + 9s=v_c]\\
      \vspace{0.2cm}
      x=9t+[9(-q-r+s)+(v_a + v_b)=v_c]\\
      \vspace{0.2cm}
      x=9(-q-r+s+t)+(v_a + v_b) = v_c \\
      \rightarrow x \equiv v_a + v_b=v_c \hspace{0.2cm} (mod\hspace{0.1cm} 9) \therefore x \in ( \phi \odot \phi)^{a+b=c}. \hspace{0.2cm} \square
     \end{matrix} 
 \end{align*}
 
\section{\texorpdfstring{Group Macro Primal ($\zeta_\phi, \circledcirc$)}{}}
 
 We define the group Macro Primal is not empty and have the next properties:
 
 \subsection{Proof:}
 Let Macro Primal be:
 \begin{align*}
     \left\{\zeta_\phi\right\}=\left\{ \Phi^Z \mid Z \in D \right\} \wedge \left\{ 0 \right\}
 \end{align*}
 
 \subsubsection{Is closed for operation:}
 
 \begin{align*}
     \phi^a \circledcirc \phi^b = ( \phi\circledcirc \phi)^{(a \circledcirc b)}
 \end{align*}
  Let's suppose
  \begin{align*}
      ( \phi\circledcirc \phi)^{(a \circledcirc b)}=( \phi' \circledcirc \phi')^{(a' \circledcirc b')}
  \end{align*}
  Then if
  
  \begin{align*}
      \begin{matrix}
      \vspace{0.2cm}
      Z_{\phi a},Z'_{\phi a} \in \phi^a \rightarrow Z_{\phi a} \equiv Z'_{\phi a}~ (mod\hspace{0.1cm} 9) \\
      \vspace{0.2cm}
      Z_{\phi b},Z'_{\phi b} \in \phi^b \rightarrow Z_{\phi b} \equiv Z'_{\phi b}~ (mod\hspace{0.1cm} 9)\\
      \vspace{0.2cm}
      \therefore \\
      \vspace{0.2cm}
      9 \mid Z_{\phi a} - Z'_{\phi a} \rightarrow Z_{\phi a} - Z'_{\phi a} = 9q / q \in \mathbb Z \\
      \vspace{0.2cm}
      9 \mid Z_{\phi b} - Z'_{\phi b}-Z_{\phi a} - Z'_{\phi b} \rightarrow Z_{\phi b} - Z'_{\phi b} = 9s / s \in \mathbb Z
      \end{matrix}
  \end{align*}
  By definition we have to:
  \begin{align*}
      ( \phi\circledcirc \phi)^{(a \circledcirc b)}=\left\{v \in \mathbb Z \mid (Z_{\phi a} \circledcirc Z_{\phi b}) ~(mod\hspace{0.1cm} 9)  \right\}
  \end{align*}
  Then we verified that
  \begin{align*}
      \begin{matrix}
       \vspace{0.2cm}
      9 \mid v - (Z_{\phi a} \circledcirc Z_{\phi b})  \\
      v-(Z_{\phi a} \circledcirc Z_{\phi b})= 9t / t \in \mathbb Z
      \end{matrix}
  \end{align*}
  then
  \begin{align*}
      \lambda) v=9t +(Z_{\phi a} \circledcirc Z_{\phi b})
  \end{align*}
  Starting from $\lambda$ we will continue the proof for the addition operation and later we will continue from it for multiplication.
  \begin{align*}
      \begin{matrix}
      \vspace{0.2cm}
       v=9t+(Z_{\phi a} + Z_{\phi b})\\
       \vspace{0.2cm}
       v=9t+(9q+ Z'_{\phi a} + 9s + Z'_{\phi b})\\
       v=9(t+q+s)+(Z'_{\phi a} + Z'_{\phi b})
      \end{matrix}
  \end{align*}
  How $(t+q+s) \in \mathbb Z $, we finally get that:
  \begin{align*}
      \begin{matrix}
      \vspace{0.2cm}
      v=9(t+q+s)+(Z'_{\phi a} + Z'_{\phi b}) \rightarrow v \equiv (Z'_{\phi a} + Z'_{\phi b})~(mod~9)\\
      \vspace{0.2cm}
      \therefore\\
      v \in (Z'_{\phi a} + Z'_{\phi b})
      \end{matrix}
  \end{align*}
  Now we continue the proof from $\lambda$ for the multiplication:
  
  \begin{align*}
      \begin{matrix}
      \vspace{0.2cm}
      v=9t+(Z'_{\phi a} * Z'_{\phi b})\\
      \vspace{0.2cm}
      v=9t+[(9q+ Z'_{\phi a})(9s+Z'_{\phi b})]\\
      \vspace{0.2cm}
      v=9t+[(81qs+9qZ'_{\phi b} +9sZ'_{\phi a} + Z'_{\phi a} * Z'_{\phi b})]\\
      \vspace{0.2cm}
      v=(9t+81qs+9qZ'_{\phi b}+9sZ'_{\phi a})+(Z'_{\phi a} * Z'_{\phi b})\\
      \vspace{0.2cm}
      v=9(t+9qs+9qZ'_{\phi b}+9sZ'_{\phi a}+(Z'_{\phi a} * Z'_{\phi b})
      \end{matrix}
  \end{align*}
  
  How $(t+9qs+9qZ'_{\phi b} + 9sZ'_{\phi a}) \in \mathbb Z $ then
  \begin{align*}
      \begin{matrix}
      \vspace{0.2cm}
      v=9(t+9qs+9qZ'_{\phi b} + 9sZ'_{\phi a})+(t+9qs+9qZ'_{\phi b} + 9sZ'_{\phi a}) \rightarrow v \equiv (Z'_{\phi a} * Z'_{\phi b})~(mod~9)\\
     \vspace{0.2cm}
      \therefore\\
      \vspace{0.2cm}
      v \in (Z'_{\phi a} * Z'_{\phi b})
      \end{matrix}
  \end{align*}
  \subsubsection{Associative property for the sum}
  
  Let
  
  \begin{align*}
      \phi^a, \phi^b, \phi^c \in \zeta_\phi \rightarrow [\phi^a+\phi^b]+\phi^c=\phi^a+[\phi^b+\phi^c]
  \end{align*}
  
  Let's supposed 
  \begin{align*}
      \begin{matrix}
      Z_{\phi a} \in \phi^a, Z_{\phi b} \in \phi^b, Z_{\phi c} \in \phi^c
      \end{matrix}
  \end{align*}
  Then
  \begin{align*}
      \begin{matrix}
      \vspace{0.2cm}
      [\phi^a+\phi^b]+\phi^c=(Z_{\phi a}+Z_{\phi b})+ \phi^c\\
      \vspace{0.2cm}
      =[(Z_{\phi a}+Z_{\phi b}+Z_{\phi c})^{a+b+c}]\\
      \vspace{0.2cm}
      =[\phi^a+(Z_{\phi b}+Z_{\phi c})^{b+c}]\\
      \vspace{0.2cm}
      =\phi^a+[Z_{\phi b}+Z_{\phi c}]^{b+c}\\
      \vspace{0.2cm}
      =\phi^a+[\phi^b+\phi^c]
      \end{matrix}
  \end{align*}
  
  \subsection{Property of existence of a neutral element}
 
  \begin{center}
      $\phi^{\pm9} \in \left\{\zeta_\phi\right\} $ is the neutral element.
  \end{center}
  
  \subsubsection{Proof:}
  \begin{align*}
        \phi^a+\phi^{\pm 9}=\left\{v \in \mathbb Z \mid v \equiv (Z_{\phi a} +9)~ (mod\hspace{0.1cm} 9)\right\}
  \end{align*}
  Then
  \begin{align*}
      \begin{matrix}
      \vspace{0.2cm}
      v-(Z_{\phi a}+9)=9t, t \in \mathbb Z\\
      \vspace{0.2cm}
      v=9t+(Z_{\phi a}+9)\\
      \vspace{0.2cm}
      v=9(t+1)+Z_{\phi a} \rightarrow v \equiv Z_{\phi a}~ (mod\hspace{0.1cm} 9) \therefore v \in (\phi \circledcirc \phi)^{a+9}
      \end{matrix}
  \end{align*}
  
  \section{\texorpdfstring{Group $\left\{\zeta_\phi\right\}+$}{}}
  Proposition of elements
  \subsection{Opposite element}
  \begin{align*}
      \phi^{-a}\in \left\{\zeta_\phi\right\}
  \end{align*}
  is the opposite element
  \subsubsection{Proof:}
  
  \begin{align*}
      (-Z_{\phi a}) \in \phi^{(-a)}
  \end{align*}
  Then by definition we can affirm that:
  
 \begin{align*}
     \begin{matrix}
     \vspace{0.2cm}
     \phi^a+\phi^{(-a)}=\left\{v \in \mathbb Z \mid v \equiv [Z_{\phi a} + (-Z_{\phi a})~(mod\hspace{0.1cm} 9)]\right\}\\
     \vspace{0.2cm}
     v=9t+[Z_{\phi a}+(-Z_{\phi a})], t \in \mathbb Z\\
     \vspace{0.2cm}
     v=9t+0 \rightarrow v \equiv 0 (mod\hspace{0.1cm} 9) \rightarrow v \equiv 9 ~(mod\hspace{0.1cm} 9)
     \end{matrix}
 \end{align*}
  Ergo
  \begin{align*}
      v \in \phi^9
  \end{align*}
  Finally
  \begin{align*}
      \phi^a+\phi^{(-a)}=\phi^{\pm9}
  \end{align*}
  
  \subsection{Cancellation law of the sum}
  Let
  \begin{align*}
      \phi^a, \phi^b, \phi^c \in \zeta_\phi~ /~\phi^a+\phi^c = \phi^b+\phi^c \rightarrow \phi^a = \phi^b
  \end{align*}
  \subsubsection{Proof:}
  As we know by definition:
  \begin{align*}
      (\phi+\phi)^{a+c}=(\phi+\phi)^{b+c}=\left\{v \in \mathbb Z \mid v \equiv Z_{\phi a} + Z_{\phi c} = Z_{\phi b} +Z_{\phi c} ~(mod \hspace{0.1cm} 9)\right\}
  \end{align*}
  From where:
  \begin{align*}
      \begin{matrix}
      \vspace{0.2cm}
      v-[(Z_{\phi a}+Z_{\phi c})=(Z_{\phi b}+Z_{\phi c})]=9t, t \in \mathbb Z\\
      \vspace{0.2cm}
      v=9t+[(Z_{\phi a}+Z_{\phi c})=(Z_{\phi b}+Z_{\phi c})]\\
      v=9t+[(Z_{\phi a}+Z_{\phi c})+(-Z_{\phi c})=(Z_{\phi b}+Z_{\phi c})+(-Z_{\phi c})]\\
      \vspace{0.2cm}
      v=9t+[(Z_{\phi a} +Z_{\phi c}-Z_{\phi c})=(Z_{\phi b}+Z_{\phi c}-Z_{\phi c})]\\
      \vspace{0.2cm}
      v= 9t+[Z_{\phi a}=Z_{\phi b}]\rightarrow v \equiv Z_{\phi a} =Z_{\phi b} ~(mod\hspace{0.1cm} 9) \therefore v \in \phi^{(a=b)}
      \end{matrix}
  \end{align*}
 
 \section{Conjecture Primal}
 \label{11}
  As we mentioned before, the Primal conjecture is a corollary of the "Acuña Theorem" in chapter \ref{acuña}, and goes in this way:
  \begin{center}
      \begin{align*}
          \phi^a \odot\phi^b \equiv\phi^c~(mod~9)~\Rightarrow~\exists~Z_a,~Z_b~/~Z_a \odot Z_b=Z_c
      \end{align*}
  \end{center}
This conjecture says, that any operations have a class or primal number as result, that means we have a solution that is congruent mod 9, but ¿if the congruence exist that means all member of the class are solutions?, or ¿One or more numbers are solution to that problem?.\\
In the case that the resulting class is congruent but not all its components are results, we can deduce that at least there are numbers within these classes that satisfy all the obtained congruence.

\begin{align*}
     \phi^a \odot\phi^b \equiv\phi^c~(mod~9)~\Rightarrow~\neg \forall~Z_a,~Z_b~/~Z_a \odot Z_b=Z_c
\end{align*}

  \section{Arithmetic Primal}
  All the early theorems are the base for the arithmetic tables presented in the paper, Acuña and Marrero (2021)\cite{tarazona2021matrize_phi}
  So we gonna show again and add some properties and new tables.
  \subsection{Primal additive Table}
  \begin{align*}
      \forall ~\phi^a+\phi^b=\phi^c
  \end{align*}
  
  \begin{center}
\begin{tabular}{|c|c|c|c|c|c|c|c|c|c|}
  \hline
        + &  $\phi^1$ & $\phi^2$ & $\phi^3$ & $\phi^4$ & $\phi^5$ & $\phi^6$ & $\phi^7$ & $\phi^8$ & $\phi^9$ \\\hline
$\phi^1$ & $2_{\phi}$ & $3_{\phi}$ & $4_{\phi}$ & $5_{\phi}$ &
$6_{\phi}$ & $7_{\phi}$ & $8_{\phi}$ & $9_{\phi}$ &
$1_{\phi}$\\
\hline

$\phi^2$ & $3_{\phi}$ & $4_{\phi}$ & $5_{\phi}$ & $6_{\phi}$ &
$7_{\phi}$ & $8_{\phi}$ & $9_{\phi}$ & $1_{\phi}$ & $2_{\phi}$
\\
\hline

$\phi^3$ & $4_{\phi}$ & $5_{\phi}$ & $6_{\phi}$ & $7_{\phi}$ &
$8_{\phi}$ & $9_{\phi}$ & $1_{\phi}$ & $2_{\phi}$ & $3_{\phi}$
\\
\hline

$\phi^4$ & $5_{\phi}$ & $6_{\phi}$ & $7_{\phi}$ & $8_{\phi}$ &
$9_{\phi}$ & $1_{\phi}$ & $2_{\phi}$ & $3_{\phi}$ & $4_{\phi}$
\\
\hline

$\phi^5$ & $6_{\phi}$ & $7_{\phi}$ & $8_{\phi}$ & $9_{\phi}$ &
$1_{\phi}$ & $2_{\phi}$ & $3_{\phi}$ & $4_{\phi}$ & $5_{\phi}$
\\
\hline

$\phi^6$ & $7_{\phi}$ & $8_{\phi}$ & $9_{\phi}$ & $1_{\phi}$ &
$2_{\phi}$ & $3_{\phi}$ & $4_{\phi}$ & $5_{\phi}$ & $6_{\phi}$
\\
\hline

$\phi^7$ & $8_{\phi}$ & $9_{\phi}$ & $1_{\phi}$ & $2_{\phi}$ &
$3_{\phi}$ & $4_{\phi}$ & $5_{\phi}$ & $6_{\phi}$ & $7_{\phi}$
\\
\hline

$\phi^8$ & $9_{\phi}$ & $1_{\phi}$ & $2_{\phi}$ & $3_{\phi}$ &
$4_{\phi}$ & $5_{\phi}$ & $6_{\phi}$ & $7_{\phi}$ & $8_{\phi}$
\\
\hline

$\phi^9$ & $1_{\phi}$ & $2_{\phi}$ & $3_{\phi}$ & $4_{\phi}$ &
$5_{\phi}$ & $6_{\phi}$ & $7_{\phi}$ & $8_{\phi}$ & $9_{\phi}$
\\\hline
\end{tabular}\\
Table 4.
\end{center}

In this table, the cycle of infinite repetitions of the primal values can be clearly appreciated.

\subsection{Primal subtraction table}

\begin{align*}
      \forall~ \phi^a-\phi^b=\phi^c
  \end{align*}

Like the subtraction is not commutative we have two tables, the first to positive results and the second for the negative results.
\begin{center}
\begin{tabular}{|c|c|c|c|c|c|c|c|c|c|}
  \hline
        - &  $\phi^1$ & $\phi^2$ & $\phi^3$ & $\phi^4$ & $\phi^5$ & $\phi^6$ & $\phi^7$ & $\phi^8$ & $\phi^9$ \\\hline
$\phi^{-1}$ & $9_{\phi}$ & $1_{\phi}$ & $2_{\phi}$ & $3_{\phi}$ &
$4_{\phi}$ & $5_{\phi}$ & $6_{\phi}$ & $7_{\phi}$ &
$8_{\phi}$\\
\hline

$\phi^{-2}$ & $8_{\phi}$ & $9_{\phi}$ & $1 _{\phi}$ & $2_{\phi}$ &
$3_{\phi}$ & $4_{\phi}$ & $5_{\phi}$ & $6_{\phi}$ & $7_{\phi}$
\\
\hline

$\phi^{-3}$ & $7_{\phi}$ & $8_{\phi}$ & $9_{\phi}$ & $1_{\phi}$ &
$2_{\phi}$ & $3_{\phi}$ & $4_{\phi}$ & $5_{\phi}$ & $6_{\phi}$
\\
\hline

$\phi^{-4}$ & $6_{\phi}$ & $7_{\phi}$ & $8_{\phi}$ & $9_{\phi}$ &
$1_{\phi}$ & $2_{\phi}$ & $3_{\phi}$ & $4_{\phi}$ & $5_{\phi}$
\\
\hline

$\phi^{-5}$ & $5_{\phi}$ & $6_{\phi}$ & $7_{\phi}$ & $8_{\phi}$ &
$9_{\phi}$ & $1_{\phi}$ & $2_{\phi}$ & $3_{\phi}$ & $4_{\phi}$
\\
\hline

$\phi^{-6}$ & $4_{\phi}$ & $5_{\phi}$ & $6_{\phi}$ & $7_{\phi}$ &
$8_{\phi}$ & $9_{\phi}$ & $1_{\phi}$ & $2_{\phi}$ & $3_{\phi}$
\\
\hline

$\phi^{-7}$ & $3_{\phi}$ & $4_{\phi}$ & $5_{\phi}$ & $6_{\phi}$ &
$7_{\phi}$ & $8_{\phi}$ & $9_{\phi}$ & $1_{\phi}$ & $2_{\phi}$
\\
\hline

$\phi^{-8}$ & $2_{\phi}$ & $3_{\phi}$ & $4_{\phi}$ & $5_{\phi}$ &
$6_{\phi}$ & $7_{\phi}$ & $8_{\phi}$ & $9_{\phi}$ & $1_{\phi}$
\\
\hline

$\phi^{-9}$ & $1_{\phi}$ & $2_{\phi}$ & $3_{\phi}$ & $4_{\phi}$ &
$5_{\phi}$ & $6_{\phi}$ & $7_{\phi}$ & $8_{\phi}$ & $9_{\phi}$
\\\hline
\end{tabular}\\
Table 5.
\end{center}

\hspace{2cm}
\begin{center}
\begin{tabular}{|c|c|c|c|c|c|c|c|c|c|}
  \hline
        + &  $\phi^{-1}$ & $\phi^{-2}$ & $\phi^{-3}$ & $\phi^{-4}$ & $\phi^{-5}$ & $\phi^{-6}$ & $\phi^{-7}$ & $\phi^{-8}$ & $\phi^{-9}$ \\\hline
$\phi^{1}$ & $-9_{\phi}$ & $-8_{\phi}$ & $-7_{\phi}$ & $-6_{\phi}$ &
$-5_{\phi}$ & $-4_{\phi}$ & $-3_{\phi}$ & $-2_{\phi}$ &
$-1_{\phi}$\\
\hline

$\phi^{2}$ & $-1_{\phi}$ & $-9_{\phi}$ & $-8 _{\phi}$ & $-7_{\phi}$ &
$-6_{\phi}$ & $-5_{\phi}$ & $-4_{\phi}$ & $-3_{\phi}$ & $-2_{\phi}$
\\
\hline

$\phi^{3}$ & $-2_{\phi}$ & $-1_{\phi}$ & $-9_{\phi}$ & $-8_{\phi}$ &
$-7_{\phi}$ & $-6_{\phi}$ & $-5_{\phi}$ & $-4_{\phi}$ & $-3_{\phi}$
\\
\hline

$\phi^{4}$ & $-3_{\phi}$ & $-2_{\phi}$ & $-1_{\phi}$ & $-9_{\phi}$ &
$-8_{\phi}$ & $-7_{\phi}$ & $-6_{\phi}$ & $-5_{\phi}$ & $-4_{\phi}$
\\
\hline

$\phi^{5}$ & $-4_{\phi}$ & $-3_{\phi}$ & $-2_{\phi}$ & $-1_{\phi}$ &
$-9_{\phi}$ & $-8_{\phi}$ & $-7_{\phi}$ & $-6_{\phi}$ & $-5_{\phi}$
\\
\hline

$\phi^{6}$ & $-5_{\phi}$ & $-4_{\phi}$ & $-3_{\phi}$ & $-2_{\phi}$ &
$-1_{\phi}$ & $-9_{\phi}$ & $-8_{\phi}$ & $-7_{\phi}$ & $-6_{\phi}$
\\
\hline

$\phi^{7}$ & $-6_{\phi}$ & $-5_{\phi}$ & $-4_{\phi}$ & $-3_{\phi}$ &
$-2_{\phi}$ & $-1_{\phi}$ & $-9_{\phi}$ & $-8_{\phi}$ & $-7_{\phi}$
\\
\hline

$\phi^{8}$ & $-7_{\phi}$ & $-6_{\phi}$ & $-5_{\phi}$ & $-4_{\phi}$ &
$-3_{\phi}$ & $-2_{\phi}$ & $-1_{\phi}$ & $-9_{\phi}$ & $-8_{\phi}$
\\
\hline

$\phi^{9}$ & $-8_{\phi}$ & $-7_{\phi}$ & $-6_{\phi}$ & $-5_{\phi}$ &
$-4_{\phi}$ & $-3_{\phi}$ & $-2_{\phi}$ & $-1_{\phi}$ & $-9_{\phi}$
\\\hline
\end{tabular}\\
Table 6.
\end{center}

Also we can obtain 0 when we subtract the same class, if and only if the number they are representing is the same.
\begin{align*}
    \begin{matrix}
    \phi^a-\phi^b = 0 \equiv v_a - v_b=0 \mid v_a=a_b \rightarrow \phi^a=\phi^b
    \end{matrix}
\end{align*}\\
These two tables in particular are complex and less intuitive, since it is an operation that depends on the size of the number and the signs to determine if the result will be positive or negative, but it also affects the result in a modular way.

\subsection{Primal Division Table}

\begin{center}
\begin{tabular}{|c|c|c|c|c|c|c|c|c|c|}
  \hline
        $\div$ &  $\phi^1$ & $\phi^2$ & $\phi^3$ & $\phi^4$ & $\phi^5$ & $\phi^6$ & $\phi^7$ & $\phi^8$ & $\phi^9$ \\\hline
$\phi^1$ & $1_{\phi}$ & $2_{\phi}$ & $3_{\phi}$ & $4_{\phi}$ &
$5_{\phi}$ & $6_{\phi}$ & $7_{\phi}$ & $8_{\phi}$ &
$9_{\phi}$\\
\hline

$\phi^2$ & $5_{\phi}$ & $1_{\phi}$ & $6_{\phi}$ & $2_{\phi}$ &
$7_{\phi}$ & $3_{\phi}$ & $8_{\phi}$ & $7_{\phi}$ & $9_{\phi}$
\\
\hline

$\phi^3$ & $\varnothing$ & $\varnothing$ & $4_{\phi};1_{\phi};7_{\phi}$ & $\varnothing$ &
$\varnothing$ & $5_{\phi};2_{\phi};8_{\phi}$ & $\varnothing$ & $\varnothing$ & $3_{\phi};9_{\phi};6_{\phi}$
\\
\hline

$\phi^4$ & $7_{\phi}$ & $5_{\phi}$ & $3_{\phi}$ & $1_{\phi}$ &
$8_{\phi}$ & $6_{\phi}$ & $4_{\phi}$ & $2_{\phi}$ & $9_{\phi}$
\\
\hline

$\phi^5$ & $2_{\phi}$ & $4_{\phi}$ & $6_{\phi}$ & $8_{\phi}$ &
$1_{\phi}$ & $3_{\phi}$ & $5_{\phi}$ & $7_{\phi}$ & $9_{\phi}$
\\
\hline

$\phi^6$ & $\varnothing$ & $\varnothing$ & $5_{\phi};2_{\phi};8_{\phi}$ & $\varnothing$ &
$\varnothing$ & $4_{\phi};1_{\phi};7_{\phi}$ & $\varnothing$ & $\varnothing$ & $3_{\phi};9_{\phi};6_{\phi}$
\\
\hline

$\phi^7$ & $4_{\phi}$ & $8_{\phi}$ & $3_{\phi}$ & $7_{\phi}$ &
$2_{\phi}$ & $6_{\phi}$ & $1_{\phi}$ & $5_{\phi}$ & $9_{\phi}$
\\
\hline

$\phi^8$ & $8_{\phi}$ & $7_{\phi}$ & $6_{\phi}$ & $5_{\phi}$ &
$4_{\phi}$ & $3_{\phi}$ & $2_{\phi}$ & $1_{\phi}$ & $9_{\phi}$
\\
\hline

$\phi^9$ & $\varnothing$ & $\varnothing$ & $\varnothing$ & $\varnothing$ &
$\varnothing$ & $\varnothing$ & $\varnothing$ & $\varnothing$ & $Z_{\phi}$
\\\hline
\end{tabular}\\
Table 7.
\end{center}

Also we know by the Acuña's Theorem, That:
\begin{align*}
    \forall\hspace{0.1cm}~\frac{\phi^a}{\phi^b}=\phi^c \rightarrow \frac{\phi^a}{\phi^c}=\phi^b \rightarrow \phi^b*\phi^c=\phi^a
\end{align*}

\subsection{Primal multiplication table}

\begin{align*}
      \forall ~\phi^a*\phi^b=\phi^c
  \end{align*}

\begin{center}
\begin{tabular}{|c|c|c|c|c|c|c|c|c|c|}
  \hline
        * &  $\phi^1$ & $\phi^2$ & $\phi^3$ & $\phi^4$ & $\phi^5$ & $\phi^6$ & $\phi^7$ & $\phi^8$ & $\phi^9$ \\\hline
$\phi^1$ & $1_{\phi}$ & $2_{\phi}$ & $3_{\phi}$ & $4_{\phi}$ &
$5_{\phi}$ & $6_{\phi}$ & $7_{\phi}$ & $8_{\phi}$ &
$9_{\phi}$\\
\hline

$\phi^2$ & $2_{\phi}$ & $4_{\phi}$ & $6_{\phi}$ & $8_{\phi}$ &
$1_{\phi}$ & $3_{\phi}$ & $5_{\phi}$ & $7_{\phi}$ & $9_{\phi}$
\\
\hline

$\phi^3$ & $3_{\phi}$ & $6_{\phi}$ & $9_{\phi}$ & $3_{\phi}$ &
$6_{\phi}$ & $9_{\phi}$ & $3_{\phi}$ & $6_{\phi}$ & $9_{\phi}$
\\
\hline

$\phi^4$ & $4_{\phi}$ & $8_{\phi}$ & $3_{\phi}$ & $7_{\phi}$ &
$2_{\phi}$ & $6_{\phi}$ & $1_{\phi}$ & $5_{\phi}$ & $9_{\phi}$
\\
\hline

$\phi^5$ & $5_{\phi}$ & $1_{\phi}$ & $6_{\phi}$ & $2_{\phi}$ &
$7_{\phi}$ & $3_{\phi}$ & $8_{\phi}$ & $4_{\phi}$ & $9_{\phi}$
\\
\hline

$\phi^6$ & $6_{\phi}$ & $3_{\phi}$ & $9_{\phi}$ & $6_{\phi}$ &
$3_{\phi}$ & $9_{\phi}$ & $6_{\phi}$ & $3_{\phi}$ & $9_{\phi}$
\\
\hline

$\phi^7$ & $7_{\phi}$ & $5_{\phi}$ & $3_{\phi}$ & $1_{\phi}$ &
$8_{\phi}$ & $6_{\phi}$ & $4_{\phi}$ & $2_{\phi}$ & $9_{\phi}$
\\
\hline

$\phi^8$ & $8_{\phi}$ & $7_{\phi}$ & $6_{\phi}$ & $5_{\phi}$ &
$4_{\phi}$ & $3_{\phi}$ & $2_{\phi}$ & $1_{\phi}$ & $9_{\phi}$
\\
\hline

$\phi^9$ & $9_{\phi}$ & $9_{\phi}$ & $9_{\phi}$ & $9_{\phi}$ &
$9_{\phi}$ & $9_{\phi}$ & $9_{\phi}$ & $9_{\phi}$ & $9_{\phi}$
\\\hline
\end{tabular}\\
Table 8.
\end{center}
\subsection{Primal power table}

\begin{center}
\begin{tabular}{|c|c|c|c|c|c|c|c|c|c|c|c|c|c|c|}
  \hline
 $()^{*}$ & $()^{2}$ & $()^{3}$ & $()^{4}$ & $()^{5}$ & $()^{6}$ & $()^{7}$ & $()^{8}$ & $()^{9}$ & $()^{10}$ & $()^{11}$ & $()^{12}$ & $()^{13}$ & $()^{14}$ & $()^{15}$ \\
\hline

  $\phi^1$ & $1_{\phi}$ & $1_{\phi}$ & $1_{\phi}$ & $1_{\phi}$ & $1_{\phi}$ & $1_{\phi}$ & $1_{\phi}$ & $1_{\phi}$ & $1_{\phi}$ & $1_{\phi}$ & $1_{\phi}$ & $1_{\phi}$ & $1_{\phi}$ & $1_{\phi}$ \\
\hline

  $\phi^2$ & $4_{\phi}$ & $8_{\phi}$ & $7_{\phi}$ & $5_{\phi}$ & $1_{\phi}$ & $2_{\phi}$ & $4_{\phi}$ & $8_{\phi}$ & $7_{\phi}$ & $5_{\phi}$ & $1_{\phi}$ & $2_{\phi}$ & $4_{\phi}$ & $8_{\phi}$ \\
\hline

  $\phi^3$ & $9_{\phi}$ & $9_{\phi}$ & $9_{\phi}$ & $9_{\phi}$ & $9_{\phi}$ & $9_{\phi}$ & $9_{\phi}$ & $9_{\phi}$ & $9_{\phi}$ & $9_{\phi}$ & $9_{\phi}$ & $9_{\phi}$ & $9_{\phi}$ & $9_{\phi}$ \\
\hline

  $\phi^4$ & $7_{\phi}$ & $1_{\phi}$ & $4_{\phi}$ & $7_{\phi}$ & $1_{\phi}$ & $4_{\phi}$ & $7_{\phi}$ & $1_{\phi}$ & $4_{\phi}$ & $7_{\phi}$ & $1_{\phi}$ & $4_{\phi}$ & $7_{\phi}$ & $1_{\phi}$ \\
\hline

  $\phi^5$ & $7_{\phi}$ & $8_{\phi}$ & $4_{\phi}$ & $2_{\phi}$ & $1_{\phi}$ & $5_{\phi}$ & $7_{\phi}$ & $8_{\phi}$ & $4_{\phi}$ & $2_{\phi}$ & $1_{\phi}$ & $5_{\phi}$ & $7_{\phi}$ & $8_{\phi}$ \\
\hline

  $\phi^6$ & $9_{\phi}$ & $9_{\phi}$ & $9_{\phi}$ & $9_{\phi}$ & $9_{\phi}$ & $9_{\phi}$ & $9_{\phi}$ & $9_{\phi}$ & $9_{\phi}$ & $9_{\phi}$ & $9_{\phi}$ & $9{\phi}$ & $9_{\phi}$ & $9_{\phi}$ \\
\hline

  $\phi^7$ & $4_{\phi}$ & $1_{\phi}$ & $7_{\phi}$ & $4_{\phi}$ & $1_{\phi}$ & $7_{\phi}$ & $4_{\phi}$ & $1_{\phi}$ & $7_{\phi}$ & $4_{\phi}$ & $1_{\phi}$ & $7_{\phi}$ & $4_{\phi}$ & $1_{\phi}$ \\
\hline

  $\phi^8$ & $1_{\phi}$ & $8_{\phi}$ & $1_{\phi}$ & $8{\phi}$ & $1_{\phi}$ & $8_{\phi}$ & $1_{\phi}$ & $8_{\phi}$ & $1_{\phi}$ & $8_{\phi}$ & $1_{\phi}$ & $8_{\phi}$ & $1_{\phi}$ & $8_{\phi}$ \\
\hline

  $\phi^9$ & $9_{\phi}$ & $9_{\phi}$ & $9_{\phi}$ & $9_{\phi}$ & $9_{\phi}$ & $9_{\phi}$ & $9_{\phi}$ & $9_{\phi}$ & $9_{\phi}$ & $9_{\phi}$ & $9_{\phi}$ & $9_{\phi}$ & $9_{\phi}$ & $9_{\phi}$ \\
  \hline
\end{tabular}\\
Table 9.
\end{center}

In this table we can see all the results for the powers from the square to the power 15, looks like a cycle appears when we reach $(\phi^Z)^8$ and reboot in  $(\phi^Z)^{14}$, so this is a open problem.

\subsubsection{Proposal Conjecture 2:"infinite cycle of powers"}
 \begin{align*}
 \forall\hspace{0.1cm} (\phi^Z)^n \hspace{0.1cm} \exists \hspace{0.1cm} (\phi^Z)^{n+6} \mid (\phi^Z)^n  \equiv (\phi^Z)^{n+6}
 \end{align*}
 This means that there is also an infinite cycle for the powers, which is symmetric and constant, this cycle would demonstrate the existence of congruence in the powers, allowing us to group polynomial equations so that they have the same modular congruence or in this case the same primal class.

 \section{Corollary of Problems}
 All this problem are in (mod\hspace{0.1cm} 9)
 \begin{itemize}
     \item \begin{align*}
        ¿ X^3+Y^3+Z^3= K \equiv X^9+Y^9+Z^9= K?
     \end{align*}
     \item \begin{align*}
         ¿(Z_{\phi a})^2+(Z_{\phi b})^3= k~for~any~specific~k
     \end{align*}
     \item \begin{align*}
     Z^n+X^n= K^n \equiv Z^{n+6}+X^{n+6}= K
 \end{align*}
 \end{itemize}

 \medskip

\bibliography{Reference}

\end{document}